\documentclass[a4paper]{article}

\usepackage{amsmath}
\usepackage{amsthm}
\usepackage{amssymb}
\usepackage{setspace}
\usepackage{url}
\usepackage{hyperref}

\makeatletter
\renewcommand{\p@enumi}{}
\renewcommand{\p@enumii}{}
\makeatother

\let\ge=\geqslant
\let\le=\leqslant
\let\emptyset=\varnothing

\theoremstyle{plain}
\newtheorem{theorem}{Theorem}[section]
\newtheorem{proposition}[theorem]{Proposition}
\newtheorem{lemma}[theorem]{Lemma}
\newtheorem{corollary}[theorem]{Corollary}

\theoremstyle{definition}
\newtheorem{defn}[theorem]{Definition}
\newtheorem{conjecture}[theorem]{Conjecture}

\makeatletter
\def\intertextshort#1{%
  \ifvmode\else\\\@empty\fi
  \noalign{%
    \penalty\postdisplaypenalty \vskip0.5\belowdisplayshortskip
    \vbox{\normalbaselines
      \ifdim\linewidth=\columnwidth
      \else \parshape\@ne \@totalleftmargin \linewidth
      \fi
      \noindent#1\par}%
    \penalty\predisplaypenalty \vskip0.5\belowdisplayshortskip
  }%
}

\makeatother

\def\nomcorr{_{\hbox{}}}

\newif\ifthesis

\DeclareMathOperator{\rank}{rank}
\DeclareMathOperator{\Sym} {Sym}

\DeclareMathOperator{\seq} {seq}
\DeclareMathOperator{\wt}  {wt}
\makeatletter
\newcommand{\Wr}{\nonscript\mskip-\medmuskip\mkern5mu\mathbin
  {\operator@font Wr}\penalty900
  \mkern5mu\nonscript\mskip-\medmuskip}
\makeatother

\DeclareSymbolFont{Shuffle}{U}{shuf}{m}{n}
\DeclareFontFamily{U}{shuf}{}
\DeclareFontShape{U}{shuf}{m}{n}{
  <5> <6> <7> <8> <9> gen * shuffle
  <10> <10.95> <12> <14.4> <17.28> <20.74> <24.88> shuffle10
  }{}
\DeclareMathSymbol{\shuffle}{\mathbin}{Shuffle}{"01}  
\DeclareMathSymbol{\cshuffle}{\mathbin}{Shuffle}{"02} 

\begin{document}

\title{Permutation Group Algebras}
\author{
  Julian D. Gilbey\\
  School of Mathematical Sciences\\
  Queen Mary, University of London\\
  Mile End Road\\
  London E1 4NS, England
}
\date{December 2001}
\maketitle


\begin{abstract}
  We consider the permutation group algebra defined by Cameron and
  show that if the permutation group has no finite orbits, then no
  homogeneous element of degree one is a zero-divisor of the algebra.
  We proceed to make a conjecture which would show that the algebra is
  an integral domain if, in addition, the group is oligomorphic.  We
  go on to show that this conjecture is true in certain special cases,
  including those of the form $H\Wr S$ and $H\Wr A$, and show that in
  the oligormorphic case, the algebras corresponding to these special
  groups are polynomial algebras.  In the $H\Wr A$ case, the algebra
  is related to the shuffle algebra of free Lie algebra theory.
  We finish by considering some integer sequences which arise from
  certain of these groups.
\end{abstract}

\let\saveclearpage=\clearpage
\let\clearpage=\relax

\newcommand{\lexlt}{<_{\text{lex}}}
\newcommand{\lexle}{\le_{\text{lex}}}
\newcommand{\lexgt}{>_{\text{lex}}}
\newcommand{\lexge}{\ge_{\text{lex}}}

\newcommand{\field}{K}
\newcommand{\alphabet}{T}
\newcommand{\ZZ}{\mathbb{Z}}
\newcommand{\DD}[2]{\Delta_{#1}^{(#2)}}
\newcommand{\KA}{\ensuremath{\field\mskip-2mu\langle\alphabet\rangle}}

\renewcommand{\(}{\textup(}
\renewcommand{\)}{\textup)}


\ifthesis
\fi

\section{Introduction}
\label{sec:pgintro}

Let $G$ be a permutation group on an (infinite) set~$\Omega$.
Cameron~\cite{Cameron} defined a commutative, associative, graded
algebra~$A(G)$ which encodes information about the action of~$G$ on
finite subsets of~$\Omega$.  It is known that this algebra has zero
divisors if $G$~has any finite orbits.  The question of what happens
when $G$~has no finite orbits is the subject of several conjectures
due to Cameron~\cite{Cameron}, and we will be exploring two of them.
The first is:

\begin{conjecture}
  \label{conj:eprime}
  If $G$ has no finite orbits, then $\varepsilon$ is a prime element
  in $A(G)$.
\end{conjecture}

Here $\varepsilon$ is a certain element in the degree one component of
the algebra, defined in section~\ref{sec:algebra}.  The following
weaker conjecture would follow from this, as we explain below.

\begin{conjecture}
  \label{conj:intdomain}
  If $G$ has no finite orbits, then $A(G)$ is an integral domain.
\end{conjecture}

The first conjecture would give us insight into the following question.
If the number of orbits of~$G$ on unordered $k$-element subsets
of~$\Omega$ is~$n_k$, then for which groups does $n_k=n_{k+1}<\infty$
hold?  We will not study this question directly here; more information
can be found in \cite{Cameron} and \cite[sect.\ 3.5]{CameronBook}.

We first show that no homogeneous element of degree one in the algebra
is a zero-divisor.  Unfortunately, it is not obvious how to extend
this argument to higher degrees.  We then go on to give a conjecture
which would, if proven, yield a proof of the weaker
conjecture~\ref{conj:intdomain}, and show that it holds in two
interesting classes of permutation groups.  It also turns out in these
two cases that the algebra~$A(G)$ is a polynomial algebra, and we
determine an explicit set of polynomial generators.  It will follow
that the stronger conjecture also holds in these cases.  Although
these results do not help to answer the question raised in the
previous paragraph (as in these cases, $n_k<n_{k+1}$ for all~$k$),
they do provide further evidence to support the conjectures.

Finally, using the inverse Euler transform, Cameron~\cite{CameronJIS}
determined the number of polynomial generators of each degree which
would be needed for certain of these algebras if they were actually
polynomial algebras.  Some of these sequences appear in The On-Line
Encyclopedia of Integer Sequences~\cite{IntSeq} in the context of free
Lie algebras.  Our work gives an explanation for the sequences
observed and the connection with free Lie algebras.

\ifthesis\clearpage\fi
\section{The graded algebra of a permutation group}
\label{sec:algebra}

We now give the definition of the algebra under consideration.  Let
$G$~be a permutation group acting on~$\Omega$.  Let $\field$~be a
field of characteristic~$0$ (either $\mathbb{Q}$ or $\mathbb{C}$ will
do).  Define $V_n(G)$ to be the $\field$-vector space of all functions
from $n$-subsets of~$\Omega$ to~$\field$ which are invariant under the
natural action of~$G$ on $n$-subsets of~$\Omega$.  Define the graded
algebra
$$A(G)=\bigoplus_{n=0}^\infty V_n(G)$$
with multiplication defined by
the rule that for any $f\in V_m(G)$ and $g\in V_n(G)$, the product
$fg\in V_{m+n}(G)$ is such that for any $(m+n)$-subset
$X\subseteq\Omega$,
$$(fg)(X)=\sum_{\substack{Y\subseteq X \\ |Y|=m}}
  f(Y)g(X\setminus Y).$$
It is easy to check that, with this multiplication, $A(G)$ is a
commutative, associative, graded algebra.

If $G$ has any finite orbits, then this algebra contains
zero-divisors.  For let $X\subseteq\Omega$ be a finite orbit, $|X|=n$,
and let $f\in V_n(G)$ be the characteristic function of this set (so
$f(X)=1$ and $f(Y)=0$ for $Y\ne X$); then clearly $f^2=0$.

Considering Conjecture~\ref{conj:intdomain}, it is clear that there
are no zero-divisors in $V_0(G)$, as multiplying by an element of
$V_0(G)$ is equivalent to multiplying by an element of~$\field$.

We also note that if there is a zero-divisor in $A(G)$, so we have
$fg=0$ with $0\ne f,g\in A(G)$, then we can consider the non-zero
homogeneous components of~$f$ and~$g$ with lowest degree; say these
are $f_m$ of degree~$m$ and $g_n$ of degree~$n$ respectively.  Then
the term of degree $m+n$ in~$fg$ will be precisely~$f_mg_n$, and as
$fg=0$, we must have $f_mg_n=0$.  So we may restrict our attention to
considering homogeneous elements, and showing that for any positive
integers $m$ and~$n$, we cannot find non-zero $f\in V_m(G)$ and $g\in
V_n(G)$ with $fg=0$.

Furthermore, we will show in the next section that $V_1(G)$ contains
no zero-divisors as long as $G$~has no finite orbits, so in
particular, the element $\varepsilon\in V_1(G)$ defined by
$\varepsilon(x)=1$ for all $x\in\Omega$ is a non-zero-divisor.  So if
$f$~is a homogeneous zero-divisor of degree~$m$, with $fg=0$, and
$g$~is homogeneous of degree $n>m$, we also have
$(\varepsilon^{n-m}f)g=0$, so $\varepsilon^{n-m}f\ne0$ is a
zero-divisor of degree~$n$.  Thus, if we wish, we can restrict our
attention to showing that, for each positive integer~$n$, we cannot
find non-zero $f,g\in V_n(G)$ with $fg=0$.

Turning now to the stronger Conjecture~\ref{conj:eprime}, we see that
the second conjecture follows from this (as in~\cite{Cameron}).  For if
$fg=0$, with $f$~and~$g$ homogeneous and non-zero, and $\deg f+\deg g$
is minimal subject to this, then $\varepsilon\mid fg$, so we can
assume $\varepsilon\mid f$ by primality.  Thus $f=\varepsilon f'$, and
$\deg f'=\deg f-1$.  Thus $\varepsilon f'g=0$, which implies $f'g=0$
by the above, contrary to the minimality of $\deg f+\deg g$.

\ifthesis\clearpage\fi
\section{The degree one case}
\label{sec:deg1}

We intend to prove the following theorem.

\begin{theorem}
  \label{thm:v1(G)}
  If $G\nomcorr$ has no finite orbits, then $V_1(G)$ contains no
  zero-divisors.
\end{theorem}

In order to prove this theorem, we will make use of a technical
proposition, which is based on a theorem of Kantor~\cite{Kantor}.  We
first quote a version of Kantor's theorem, as we will have use for it
later.

\begin{proposition}
  \label{thm:kantor}
  Let $0\le e<f\le d-e$.  Let $X\nomcorr$ be a set with $|X|=d$.  We
  define $(E,F)$ for subsets $E,F\subset X\nomcorr$ with $|E|=e$ and
  $|F|=f\nomcorr$ by
  $$(E,F)=
    \begin{cases}
      1 &\text{if $E\subset F$}\\
      0 &\text{otherwise,}
    \end{cases}
  $$
  and the matrix $M=((E,F))$, where the rows of~$M\nomcorr$ are
  indexed by the $e$-subsets of~$X\nomcorr$ and the columns by the
  $f$-subsets.

  Then $\rank M=\binom{d}{e}$.
\end{proposition}

The extension of this result is as follows.

\begin{proposition}
  \label{prop:kantor-ext}
  Let $0\le e<f\le d-2e$.  Let $X\nomcorr$ be a set with $|X|=d$, and
  let $E_0\subset X\nomcorr$ with $|E_0|=e$ be a distinguished subset
  of~$X\nomcorr$.  Let $w\nomcorr$~be a weight function on the
  $(f-e)$-subsets of~$X\nomcorr$ with values in the
  field~$\field\nomcorr$, satisfying the condition that $w(X')=1$
  whenever $X'$~is an $(f-e)$-subset of~$X\nomcorr$ such that
  $X'\nsubseteq E_0$.  We define $(E,F)$ for subsets $E,F\subset
  X\nomcorr$ with $|E|=e$ and $|F|=f$ by
  $$(E,F)=
    \begin{cases}
      w(F\setminus E)&\text{if $E\subset F$}\\
      0              &\text{otherwise,}
    \end{cases}
    $$
    and the matrix $M=((E,F))$, where the rows of~$M\nomcorr$ are
    indexed by the $e$-subsets of~$X\nomcorr$ and the columns by the
    $f\nomcorr$-subsets.

  Then $\rank M=\binom{d}{e}$.
\end{proposition}

\begin{proof}[Proof of theorem \ref{thm:v1(G)}]
  Let $g\in V_1(G)$ with $g\ne0$, and assume $h\in V_n(G)$ with
  $n\ge1$ and $gh=0$ (the $n=0$ case has been dealt with in
  section~\ref{sec:algebra}).  We must show that $h=0$, so that for
  any $Y\subset\Omega$ with $|Y|=n$, we have $h(Y)=0$.  We assume that
  a set $Y$~has been fixed for the remainder of this proof.

  Since $g\ne0$, there exists some (infinite) orbit
  $\Delta\subseteq\Omega$ on which $g$~is non-zero; multiplying by a
  scalar if necessary, we may assume that $g(\delta)=1$ for all
  $\delta\in\Delta$.  Pick $X\subset\Omega$ with $|X|=3n+1$, $Y\subset
  X$ and $X\setminus Y\subset\Delta$.

  Now for any $(n+1)$-subset $F\subset X\nomcorr$, we have $(hg)(F)=0$ as
  $gh=hg=0$, so that
  $$(hg)(F)=\sum_{\substack{E\subset F \\ |E|=n}}
  h(E)g(F\setminus E)=0.$$
  This can be thought of as a system of linear equations in the
  unknowns $h(E)$ for $E\subset X\nomcorr$, $|E|=n$, with the
  matrix~$M=(m_{EF})$ given by $m_{EF}=g(F\setminus E)$ if $E\subset
  F$, and $m_{EF}=0$ otherwise.
  
  This is precisely the situation of the proposition if we let $e=n$,
  $f=n+1$ (so that $f-e=1$), $d=3n+1$, $E_0=Y$ and
  $w(\alpha)=g(\alpha)$; note that $w(\alpha)=1$ whenever
  $\alpha\notin E_0$.  (We write $g(\alpha)$ instead of the more
  correct $g(\{\alpha\})$; no confusion should arise because of this.)
  Thus $\rank M=\binom{d}{e}$ and the system of equations has a unique
  solution, which must be $h(E)=0$ for all~$E\subset X$ with $|E|=n$,
  as this is a possible solution.  In particular, this means that
  $h(Y)=0$, and since $Y$~was chosen arbitrarily, it follows that
  $h=0$.

  Hence $g$ is not a zero-divisor.
\end{proof}

\begin{proof}[Proof of proposition \ref{prop:kantor-ext}]
  \newcommand{\EEi}{\ensuremath{\mathcal{E}(E',i)}}
  Let $R(E)$ be the row of~$M$ corresponding to~$E$.  $M$~has
  $\binom{d}{e}$~rows, so we must show that the rows are linearly
  independent.  We thus assume that there is a linear dependence among
  the rows of~$M$, so
  \begin{equation}
    \label{eq:1}
    R(E^*)=\sum_{E\ne E^*} a(E)R(E)
  \end{equation}
  for some $e$-set $E^*$ and some $a(E)\in\field$.  We first note that
  $R(E^*)$ itself is non-zero: this follows as we can pick
  some~$F\supset E^*$ with $F\setminus E^*\nsubseteq E_0$; for
  this~$F$, we have $(E^*,F)=1$.

  Let $\Gamma$ be the subgroup of $\Sym(X)$ which stabilises~$E_0$
  pointwise and $E^*$~setwise.  If $\sigma\in\Gamma$, then
  $$(E^\sigma,F^\sigma)=
    \begin{cases}
    w((F\setminus E)^\sigma)=w(F\setminus E) &\text{if $E\subset F$}\\
    0                                        &\text{otherwise;}    
    \end{cases}
  $$
  either way, $(E^\sigma,F^\sigma)=(E,F)$.  (For the result
  $w((F\setminus E)^\sigma)=w(F\setminus E)$, note that both sides are
  equal to~$1$ unless $F\setminus E\subseteq E_0$, in which case
  $\sigma$~fixes this set pointwise.)  Thus \eqref{eq:1}~implies that,
  for all~$F$,
  \begin{align}
    \label{eq:2}
    (E^*,F)=(E^*,F^\sigma)
      &= \sum_{E\ne E^*} a(E^\sigma)(E^\sigma,F^\sigma)\notag\\
      &= \sum_{E\ne E^*} a(E^\sigma)(E,F).\notag\\
  \intertextshort{Thus}
    R(E^*)&=\sum_{E\ne E^*} a(E^\sigma)R(E).\notag\\
  \intertextshort{It follows that}
    |\Gamma|\,R(E^*)
      &= \sum_{\sigma\in\Gamma}\sum_{E\ne E^*} a(E^\sigma)R(E)\notag\\
      &= \sum_{E\ne E^*} R(E) \sum_{\sigma\in\Gamma} a(E^\sigma).
  \end{align}
  
  We now consider the orbits of~$\Gamma$ on the $e$-subsets
  of~$X\nomcorr$, excluding~$E^*$\negthinspace.  The $e$-sets $E_1$
  and $E_2$ will lie in the same orbit if and only if $E_1\cap
  E_0=E_2\cap E_0$ and $|E_1\cap E^*|=|E_2\cap E^*|$.  Thus every
  orbit is described by a subset $E'\subseteq E_0$ and an integer
  $0\le i\le e-1$.  (We cannot have $i=e$, as we are excluding~$E^*$
  from consideration.)  Clearly not all possible pairs $(E',i)$ will
  actually correspond to an orbit (it is not hard to see that
  necessary and sufficient conditions for this are $|E'\cap E^*|\le
  i\le\min\{e-1,\allowbreak e-|E'\setminus E^*|\}$), so that whenever
  we consider or sum over such pairs below, we implicitly restrict
  attention to those which correspond to an orbit.  In such cases, we
  write $\EEi$ for the orbit.  Also, for each such pair, pick some
  $E(E',i)\in\EEi$.  Then \eqref{eq:2}~implies
  \begin{align}
    \label{eq:3}
    |\Gamma|\,R(E^*)
      &= \sum_{(E',i)}\sum_{E\in\EEi} R(E) \sum_{\sigma\in\Gamma}
           a(E^\sigma)\notag\\
      &= \sum_{(E',i)}\sum_{E\in\EEi} R(E) \sum_{\sigma\in\Gamma}
           a(E(E',i)^\sigma)\notag\\
      &= \sum_{(E',i)}\sum_{\sigma\in\Gamma} a(E(E',i)^\sigma)
           \sum_{E\in\EEi} R(E)\notag\\
  \intertextshort{so that}
    R(E^*)&=\sum_{(E',i)} b(E',i) \sum_{E\in\EEi} R(E)
  \end{align}
  with $b(E',i)\in\field$, and clearly not all of the $b(E',i)$ can be
  zero as $R(E^*)$ is not zero.
  
  We define a total order on the pairs $(E',i)$ as follows.  Extend
  the partial order given by~$\subseteq$ on the subsets of~$E_0$ to a
  total order~$\le$, and then define $(E',i)\le(E'',j)$ if $E'<E''$ or
  $E'=E''$ and $i\le j$.  We now proceed to derive a contradiction by
  showing that \eqref{eq:3}~leads to a system of linear equations for
  the $b(E',i)$ which is triangular under this total order, with
  non-zero diagonal entries, and deduce that all of the $b(E',i)$ must
  be zero.
  
  Let $(\bar E,n)$ be a pair corresponding to an orbit.  Since
  $2e+f\le d$, there exists an $f$-set $F(\bar E,n)$ satisfying
  $F(\bar E,n)\cap E_0=\bar E$ and $|F(\bar E,n)\cap E^*|=n$.  (Simply
  take $E(\bar E,n)$ and adjoin $f-e$ points lying in
  $X\setminus(E_0\cup E^*)$.)  As $n\le e-1$, it follows that $F(\bar
  E,n)\nsupseteq E^*$, so $(E^*,F(\bar E,n))=0$.  Hence
  by~\eqref{eq:3}, we have
  \begin{equation}
    \label{eq:4}
    0=\sum_{(E',i)} b(E',i) \sum_{E\in\EEi} (E,F(\bar E,n))
  \end{equation}
  for all such pairs $(\bar E,n)$.
  
  We note that $F(\bar E,n)\cap E_0=\bar E$, and further that
  $E\in\EEi$ implies that $E\cap E_0=E'$; thus for the term $(E,F(\bar
  E,n))$ in equation~\eqref{eq:4} to be non-zero, where
  $E\in\mathcal{E}(E',i)$, we require $E'\subseteq \bar E$, hence also
  $E'\le\bar E$.  Furthermore, if $(E,F(\bar E,n))\ne0$, we must have
  $i\le n$ as $E\subset F(\bar E,n)$.  Thus if $(\bar E,n)<(E',i)$, we
  have
  \begin{equation}
    \label{eq:5}
    \sum_{E\in\EEi} (E,F(\bar E,n))=0.
  \end{equation}
  
  Also, there is an $e$-set $E\subset F(\bar E,n)$ satisfying $E\cap
  E^*=F(\bar E,n)\cap E^*$ and $E\cap E_0=F(\bar E,n)\cap E_0=\bar
  E$; just take the union of $\bar E$ with $F(\bar E,n)\cap E^*$ and
  sufficiently many remaining points of $F(\bar E,n)$.  For each
  such~$E$, we have $F(\bar E,n)\setminus E \nsubseteq E_0$, so
  $(E,F(\bar E,n))=1$.  Since $\field$~has characteristic zero, we
  deduce that
  \begin{equation}
    \label{eq:6}
    \sum_{E\in\mathcal{E}(\bar E,n)} (E,F(\bar E,n))\ne0,
  \end{equation}
  as the sum is over all sets of precisely this form.
  
  It then follows from \eqref{eq:4} and \eqref{eq:5} that for each
  pair $(\bar E,n)$:
  $$0=\sum_{(E',i)\le(\bar E,n)} b(E',i) \sum_{E\in\EEi} (E,F(\bar
  E,n)).$$
  Now this is a system of linear equations in the unknowns $b(E',i)$
  which is lower triangular.  Also, by~\eqref{eq:6}, the diagonal
  entries are non-zero.  It follows that the unique solution to this
  system is that all of the $b(E',i)$ are zero, which provides the
  required contradiction to equation~\eqref{eq:3} above.
\end{proof}

\ifthesis\clearpage\fi
\section{Oligomorphic-type cases: our conjecture}
\label{sec:conjecture}

\subsection{Ramsey orderings on orbits of $n$-sets}
\label{sec:ramsey-order}

Cameron proved the following Ramsey-type result in
\cite[Prop.\ 1.10]{CameronBook}.

\begin{lemma}
  \label{lemma:cameron}
  Suppose that the $n$-sets of an infinite set~$X\nomcorr$ are
  coloured with $r\nomcorr$~colours, all of which are used.  Then
  there is an ordering $c_1$, \dots, $c_r\!$ of the colours and
  infinite subsets $X_1$, \ldots, $X_r$, such that $X_i$~contains an
  $n$-set of colour~$c_i$ but no set of colour~$c_j$ for $j>i$.
\end{lemma}

We use this as the inspiration for the following definition.  If $G$
is a permutation group on~$\Omega$, we say that the orbits of~$G$ on
$n$-sets of~$\Omega$ can be \textit{Ramsey ordered} if, given any
finite $N>n$, there is an ordering of the orbits~$c_\alpha$,
$\alpha\in\mathcal{A}$, where $\mathcal{A}$~is a well-ordered set, and
a corresponding sequence of (possibly infinite)
subsets~$X_\alpha\subseteq\Omega$ with $|X_\alpha|\ge N$, and such
that $X_\alpha$~contains an $n$-set in the orbit~$c_\alpha$ but no
$n$-set in an orbit~$c_\beta$ for $\beta>\alpha$.  (We can take
$\mathcal{A}$~to be a set of ordinals with the $\in$-ordering if we
wish; this is the reason for using Greek letters.)  This pair of
sequences forms a \textit{Ramsey ordering}.  While the particular
Ramsey ordering may depend on~$N$, we do not usually mention~$N$
unless we have to.  The reader may think throughout of $N$~having a
very large finite value.  It turns out that this makes certain
constructions below simpler than if we required the~$X_\alpha$ to be
infinite sets.

Not every permutation group has such an ordering.  For example, in the
regular action of~$\ZZ$ on~$\ZZ$, there is no set with more than two
elements, all of whose $2$-subsets are in the same orbit, so there
cannot be a Ramsey ordering on $2$-subsets.  However, Cameron's result
implies that if $G$~is \textit{oligomorphic} (that is, there are only
finitely many orbits on $n$-sets for each~$n$), then the orbits of~$G$
on $n$-sets can be Ramsey ordered for each~$n$.

It turns out that Ramsey orderings on $n$-sets naturally yield Ramsey
orderings on $m$-sets whenever $m<n$.

\begin{proposition}
  \label{prop:induced-orbits}
  Let $G\nomcorr$ be a permutation group acting on an infinite
  set~$\Omega$.  Let $m<n$ be positive integers, and assume that the
  $n$-set orbits of~$G\nomcorr$ can be Ramsey ordered, say
  $c_\alpha\!$ and $X_\alpha\!$ with $\alpha\in\mathcal{A}$ are a
  Ramsey-ordering with $N\ge m+n$.  Then this ordering induces a
  Ramsey ordering on the $m$-set orbits as follows.  There is a
  subset\/ $\mathcal{B}\subseteq\mathcal{A}$ and a labelling of the
  $m$-set orbits as $d_\beta$, $\beta\in\mathcal{B}$, such that for
  each $\beta\in\mathcal{B}$, an $m$-set in the orbit~$d_\beta$
  appears in~$X_\beta$, and that for each $\alpha\in\mathcal{A}$,
  $X_\alpha$~contains no $m$-sets in the orbit~$d_\beta$ for
  $\beta>\alpha$.
\end{proposition}

We call the ordering of orbits $d_\beta$, $\beta\in\mathcal{B}$
together with the corresponding sets~$X_\beta$ given by this
proposition the \textit{induced Ramsey ordering}.  Note that we use
the same parameter~$N$ in both orderings.

The proof uses the following application of Kantor's theorem
(Proposition~\ref{thm:kantor} above), shown to me by Peter Cameron.

\begin{lemma}
  \label{lem:colour-types}
  Let $m<n$ be positive integers, and let $X\nomcorr$~be a finite set
  with $|X|\ge m+n$.  Let the $m$-sets of~$X\nomcorr$ be coloured with
  colours from the set\/~$\mathbb{N}$.  Given an $n$-subset
  of~$X\nomcorr$, we define its \emph{colour-type} to be the multiset
  of colours of its $\binom{n}{m}$ $m$-subsets.  Then the number of
  distinct $m$-set colours used in~$X\nomcorr$ is less than or equal
  to the number of distinct colour-types among the $n$-subsets
  of~$X\nomcorr$.
\end{lemma}

\begin{proof}
  We note that only a finite number of colours appear among the
  $m$-subsets of~$X\nomcorr$, as they are finite in number.  Without
  loss of generality, we may assume that the colours used are
  precisely $1$,~$2$, \dots,~$s$.
  
  As in Kantor's theorem (Proposition~\ref{thm:kantor}), we let $M$~be
  the incidence matrix of the $m$-subsets versus $n$-subsets of~$X$.
  By that theorem, as $m<n$ and $|X|\ge m+n$, this matrix has
  rank~$\binom{|X|}{m}$, which equals the number of rows in the
  matrix.  Thus, by the rank-nullity theorem, $M$~represents an
  injective linear transformation.
  
  Now for each $i=1$,~\dots, $s$, let $v_i$~be the row vector, with
  entries indexed by the $m$-subsets of~$X\nomcorr$, whose $j$-th
  entry is~$1$ if the $j$-th $m$-subset has colour~$i$, and $0$~if it
  does not.  Then $v_iM$ is a row vector, indexed by the $n$-subsets
  of~$X\nomcorr$, whose $k$-th entry is the number of $m$-subsets of
  the $k$-th $n$-subset which have colour~$i$.
  
  Consider now the matrix~$M'$ whose rows are $v_1M\nomcorr$, \dots,
  $v_sM\nomcorr$.  Note that the $k$-th column of this matrix gives
  the colour-type of the $k$-th $n$-subset of~$X\nomcorr$.  Its rank
  is given by
  $$\rank M'=\dim\,\langle v_1M,\dotsc,v_sM\rangle=\dim\,\langle
  v_1,\dotsc,v_s\rangle=s,$$
  as $M$~represents an injective linear transformation, and the
  $s$~vectors $v_1$,~\dots, $v_s$ are clearly linearly independent.
  Now since the row rank and column rank of a matrix are equal, we
  have $s=\rank M'\le{}$number of distinct columns in~$M'$, which is
  the number of $n$-set colour-types in~$X\nomcorr$.  Thus the number
  of $m$-set colours appearing in~$X$ is less than or equal to the
  number of $n$-set colour-types in~$X\nomcorr$, as we wanted.
\end{proof}

\begin{proof}[Proof of Proposition~\ref{prop:induced-orbits}]
  Let $c_\alpha$ be any $n$-set orbit, and let $X$~be a representative
  of this orbit.  We observe that the multiset of $m$-set orbits
  represented by the $\binom{m}{n}$ $m$-subsets of~$X$ is independent
  of the choice of~$X$ in this orbit.  (For let $\bar X$~be another
  representative of the orbit~$c_\alpha$, with $\bar X=g(X)$, where
  $g\in G$.  Then the set of $m$-subsets of~$X$ is mapped to the set
  of $m$-subsets of~$\bar X$ by~$g$, and so the multisets of $m$-set
  orbits represented by these two sets are identical.)  In particular,
  we may say that an $n$-set orbit contains an $m$-set orbit, meaning
  that any representative of the $n$-set orbit contains a
  representative of the $m$-set orbit.
  
  We first claim that every $m$-set orbit appears in some~$X_\alpha$:
  take a representative of an $m$-set orbit, say $Y\subset\Omega$.
  Adjoin a further $n-m$ elements to get an $n$-set~$\bar X\nomcorr$.
  This $n$-set lies in some orbit, so there is a representative of
  this orbit in one of the~$X_\alpha$, say $X\subset X_\alpha$.
  Then this~$X_\alpha$ contains a representative of our $m$-set orbit
  by the above argument, as we wished to show.

  Now if $Y\subset\Omega$~is a representative of an $m$-set orbit, we
  set
  $$\beta_Y=\min\,\{\,\alpha:\text{$g(Y)\subset X_\alpha$ for some
    $g\in G$}\,\}.$$
  Note that this implies that the $m$-set orbit
  containing~$Y$ is contained in $c_{\beta_Y}$ but not in $c_\alpha$
  for any $\alpha<\beta_Y$.  We set
  $\mathcal{B}=\{\,\beta_Y:Y\subset\Omega\ \text{and}\ |Y|=m\,\}$, and
  if $Y$~is an $m$-set, then we set $d_{\beta_Y}\!$~to be the orbit
  of~$Y\nomcorr$.  We claim that $\mathcal{B}$~satisfies the
  conditions of the proposition with this orbit labelling.  Certainly
  an $m$-set in the orbit~$d_{\beta_Y}$ appears in~$X_{\beta_Y}$ for
  each~$Y\nomcorr$, by construction, and for each
  $\alpha\in\mathcal{A}$, $X_\alpha$~contains no $m$-sets in the
  orbit~$d_\beta$ for $\beta>\alpha$, again by construction.  However,
  for $d_{\beta_Y}$ to~be well-defined, we require that
  $\beta_{Y_1}\ne\beta_{Y_2}$ if $Y_1$ and~$Y_2$ lie in distinct
  orbits.  We now show this to be the case by demonstrating that given
  any $\alpha_0\in\mathcal{A}$, there can only be one $m$-set orbit
  appearing in~$c_{\alpha_0}$ which has not appeared in any~$c_\alpha$
  with $\alpha<\alpha_0$.
  
  So let $\alpha_0\in\mathcal{A}$, and let $X\subseteq X_{\alpha_0}$
  have size $m+n$ and contain an $n$-set in the orbit~$c_{\alpha_0}$.
  By the observation we made above, namely that the $m$-set orbits
  appearing in an $n$-set are independent of the choice of the $n$-set
  in its $n$-set orbit, it suffices to show that our set $X$~contains
  at most one new $m$-set orbit.  To use the lemma, we colour the
  $m$-subsets of~$X$ as follows.  If $Y$~is an $m$-set with
  $\beta_Y<\alpha_0$, then $Y$~is given colour~$1$.  Those $Y\subset
  X$ with $\beta_Y=\alpha_0$ are given the colours $2$,~$3$,~\dots,
  with a distinct colour per $m$-set orbit.  (Note that any $Y\subset
  X$ has $\beta_Y\le\alpha_0$, as all $n$-subsets of~$X_{\alpha_0}$
  lie in orbits~$c_\alpha$ with $\alpha\le\alpha_0$.)
  
  We now consider the possible colour-types of the $n$-sets of~$X$.
  Note first that since the $m$-sets in a given $m$-set orbit all have
  the same colour, the colour-type of an $n$-set depends only upon the
  $n$-set orbit in which it lies.  There is some $n$-subset of~$X$ in
  the orbit~$c_{\alpha_0}$ by construction, and this has a certain
  colour-type.  Any other $n$-subset $\tilde X\subset X$ is either in
  the same orbit~$c_{\alpha_0}$, and so has the same colour-type, or
  it is in some other orbit~$c_\alpha$ with $\alpha<\alpha_0$.  In
  the latter case, every $m$-subset $Y\subset\tilde X$ must have
  $\beta_Y\le\alpha<\alpha_0$, and so it has colour~$1$.  Thus the
  colour-type of such an $n$-set must be the multiset
  $[1,1,\dotsc,1]$.
  
  If every $n$-subset of~$X$ is in the orbit~$c_{\alpha_0}$, then
  there is only one colour-type, and so there can only be one $m$-set
  colour in~$X$ by the lemma, that is, only one $m$-set orbit with
  $\beta_Y=\alpha_0$.  On the other hand, if $X$~contains an $n$-set
  in an orbit~$c_\alpha$ with $\alpha<\alpha_0$, then there are at
  most two colour-types in~$X$: the all-$1$ colour-type and the
  colour-type of $c_{\alpha_0}$.  Thus, by the lemma, $X$~contains at
  most two $m$-set colours.  Colour~$1$ appears in $c_\alpha$, and so
  there is at most one other colour present, that is, there is at most
  one $m$-set orbit with $\beta_Y=\alpha_0$.  Thus $d_{\beta_Y}$~is
  well-defined on $m$-set orbits, and we are done.
\end{proof}

\subsection{The Ramsey-ordering conjecture}
\label{sec:ramsey-conjecture}

Let $G$ be a permutation group on~$\Omega$ and let $m$ and~$n$ be
positive integers.  Let $d$~be an $m$-set orbit and $e$~an $n$-set
orbit.  If $c$~is an $(m+n)$-set orbit, then we say that
\textit{$c$~contains a $d\cup e$ decomposition} if an $(m+n)$-set~$X$
in the orbit~$c$ can be written as $X=X_m\cup X_n$ with $X_m$~in~$d$
and $X_n$~in~$e$.  We can easily show using a theorem of
P.\,M.~Neumann that if $G$~has no finite orbits, then for every pair
$(d,e)$, there exists an $(m+n)$-set orbit~$c$ containing a $d\cup e$
decomposition, as follows.

Neumann~\cite{Neumann} proved the following: Let $G$~be a permutation
group on~$\Omega$ with no finite orbits, and let $\Delta$ be a finite
subset of~$\Omega$.  Then there exists $g\in G$ with $g\Delta\cap
\Delta=\emptyset$.  It follows trivially that if $Y$ and~$Z$ are
finite subsets of~$\Omega$, then there exists $g\in G$ with $gY\cap
Z=\emptyset$ (just take $\Delta=Y\cup Z$).  In our case, let $X_m$
and~$X_n$ be representatives of~$d$ and~$e$ respectively.  Then there
exists $g\in G$ with $gX_m\cap X_n=\emptyset$, and $gX_m\cup X_n$ is
an $(m+n)$-set with the required decomposition, hence we can take
$c$~to be its orbit.

We will be considering groups~$G$ which have a Ramsey ordering on
their $(m+n)$-set orbits.  Let $c_\alpha$, $\alpha\in\mathcal{A}$ be
the ordering on $(m+n)$-sets, and let $d_\beta$, $\beta\in\mathcal{B}$
and $e_\gamma$, $\gamma\in\mathcal{C}$ be the induced Ramsey orderings
on $m$- and $n$-sets respectively (where we assume $N$~is sufficiently
large).  We then define
$$\beta\vee\gamma = \min\,\{\,\alpha:\text{$c_\alpha$~contains a
  $d_\beta\cup e_\gamma$ decomposition}\,\}.$$

Here is our main conjecture.

\begin{conjecture}
  \label{conj:main}
  Let $G$ be a permutation group on~$\Omega$ with no finite orbits and
  for which the orbits on $n$-sets can be Ramsey ordered for
  every~$n$.  Then given positive integers $m$ and~$n$, there exists
  some Ramsey ordering of the orbits on $(m+n)$-sets with $N\ge
  2(m+n)$, say $c_\alpha$, $\alpha\in\mathcal{A}$ with corresponding
  sets $X_\alpha\subseteq\Omega$, which induces Ramsey orderings
  $d_\beta$, $\beta\in\mathcal{B}$ and $e_\gamma$,
  $\gamma\in\mathcal{C}$ on the $m$-set orbits and $n$-set orbits
  respectively, and which satisfies the following conditions for all
  $\beta,\beta'\in\mathcal{B}$ and $\gamma,\gamma'\in\mathcal{C}$:
  $$\text{$\beta\vee\gamma<\beta'\vee\gamma$ if\/
  $\beta<\beta'$ \quad and \quad $\beta\vee\gamma<\beta\vee\gamma'$
  if\/ $\gamma<\gamma'$.}$$
\end{conjecture}

Note that the conditions of this conjecture also imply that if
$\beta<\beta'$ and $\gamma<\gamma'$, then
$\beta\vee\gamma<\beta\vee\gamma'<\beta'\vee\gamma'$, so that
$\beta\vee\gamma\le\beta'\vee\gamma'$ implies that either
$\beta<\beta'$ or $\gamma<\gamma'$ or
$(\beta,\gamma)=(\beta',\gamma')$.

Given this conjecture, it is easy to show that $A(G)$ is an integral
domain for such groups.  For if $fg=0$ with $0\ne f\in V_m(G)$ and
$0\ne g\in V_n(G)$, let $\beta_0$ be such that $f(d_\beta)=0$ for
$\beta<\beta_0$ but $f(d_{\beta_0})\ne0$, and let $\gamma_0$ be such
that $g(e_\gamma)=0$ for $\gamma<\gamma_0$ but $g(e_{\gamma_0})\ne0$.
(We write $f(d_\beta)$ to mean the value of $f(Y)$ where $Y$~is any
representative of the orbit~$d_\beta$, and so on.)  Letting
$\alpha_0=\beta_0\vee\gamma_0$, we can consider $fg(c_{\alpha_0})$.
Now since $fg=0$, this must be zero, but we can also determine this
explicitly.  Letting $X$~be a representative of $c_{\alpha_0}$, we
have
$$fg(c_{\alpha_0})=fg(X)=\sum_{\substack{Y\subset X\\|Y|=m}}
f(Y)g(X\setminus Y).$$
Every term in the sum is of the form
$f(d_\beta)g(e_\gamma)$ where $d_\beta\cup e_\gamma$ is a
decomposition of~$c_{\alpha_0}$, so that
$\beta\vee\gamma\le\alpha_0=\beta_0\vee\gamma_0$.  But by the
conjecture, this implies that except for terms of the form
$f(d_{\beta_0})g(e_{\gamma_0})\ne0$, every term either has
$\beta<\beta_0$ so that $f(d_\beta)=0$, or $\gamma<\gamma_0$ so that
$g(e_\gamma)=0$, and hence every one of these terms is zero.  Since
there exist terms of the form $f(d_{\beta_0})g(e_{\gamma_0})$ by the
choice of~$\alpha_0$, we must have $fg(c_{\alpha_0})\ne0$.  But this
contradicts $fg=0$, and so $A(G)$~is an integral domain.

Recall from section~\ref{sec:algebra} that we can assume $m=n$ when
showing that $A(G)$ is an integral domain (that is, $fg=0$ where
$f,g\in V_n(G)$ implies $f=0$ or $g=0$); hence we can restrict
ourselves to proving the conjecture in the case $m=n$ if this is
easier.

\ifthesis\clearpage\fi
\section{Special cases (I): Wreath-$S$-like groups}
\label{sec:special-I}

\subsection{Notational conventions}
\label{sec:pgnotation}

We gather here some notation that we will be using for the rest of
this \ifthesis part of the thesis\else paper\fi.

We will make use of the lexicographical order on finite sequences and
multisets, which we define as follows.  Let $(X,<)$~be a totally
ordered set.  If $x=(x_1,\dotsc,x_r)$ and $y=(y_1,\dots,y_s)$ are two
ordered sequences of elements of~$X\nomcorr$, then we say that $x$~is
lexicographically smaller than~$y$, written $x\lexlt y$, if there is
some~$t$ with $x_i=y_i$ for all $i<t$, but either $x_t<y_t$ or
$r+1=t\le s$.  If we now take a finite multiset of elements
of~$X\nomcorr$, say~$M$, we write~$\seq(M)$ to mean the sequence
obtained by writing the elements of~$M$ (as many times as they appear
in~$M$) in decreasing order.  Then if $M_1$ and~$M_2$ are finite
multisets, we define $M_1\lexlt M_2$ to mean
$\seq(M_1)\lexlt\seq(M_2)$.  Note that $\lexlt$~is a total order on
the set of finite multisets, for $\seq(M_1)=\seq(M_2)$ if and only if
$M_1=M_2$.  If we need to explicitly list the elements of a multiset,
we will write $[x_1,x_2,\dotsc]$.  We write $M_1+M_2$ for the multiset
sum of the multisets~$M_1$ and~$M_2$, so if $M_1=[x_1,\dotsc,x_r]$ and
$M_2=[y_1,\dotsc,y_s]$, then
$M_1+M_2=[x_1,\dotsc,x_r,y_1,\dotsc,y_s]$.

In the following sections, we will talk about a set of
\textit{connected blocks} for a permutation group, the idea being that
every orbit will correspond to a multiset or sequence of connected
blocks.  The choice of terminology will be explained below, and is not
related to blocks of imprimitivity.  Also, the individual words
``connected'' and ``block'' have no intrinsic meaning in the context
of the definitions in this \ifthesis thesis\else paper\fi.  Every
connected block has a positive integral weight (for which we write
$\wt(\Delta)$), and the weight of a sequence or multiset of connected
blocks is just the sum of weights of the individual connected blocks.
We well-order the connected blocks of each weight, and denote the
connected blocks of weight~$i$ by~$\DD ij$, where $j$~runs through
some well-ordered indexing set.  Without loss of generality, we assume
that $\DD11$~is the least connected block of weight~$1$.  We then
define a well-ordering on all connected blocks by $\DD ij<\DD{i'}{j'}$
if $i<i'$ or $i=i'$ and $j<j'$.  Using this ordering, we can then talk
about the lexicographic ordering on sequences or multisets of
connected blocks.

\subsection{Wreath-$S$-like groups}
\label{sec:wreath-S}

Our prototypical family of groups for this class of groups are those
of the form $G=H\Wr S$, where $H$~is a permutation group on~$\Delta$
and $S=\Sym(\ZZ)$, the symmetric group acting on a countably infinite
set (we take the integers for convenience).  The action is the
imprimitive one, so $G$~acts on $\Omega=\Delta\times\ZZ$.  We extract
those features of this group which are necessary for the proof below
to work.

\begin{defn}
  We say that a permutation group~$G$ on~$\Omega$ is
  \emph{wreath-$S$-like} if there is a set of connected
  blocks~$\{\DD ij\}$ and a bijection~$\phi$ from the set of orbits
  of~$G$ on finite subsets of~$\Omega$ to the set of all finite
  multisets of connected blocks, with the bijection satisfying the following
  conditions (where we again blur the distinction between orbits and
  orbit representatives):
  \begin{enumerate}
  \item  If $Y\subset\Omega$ is finite, then $\wt(\phi(Y))=|Y|$.
  \item If $Y\subset\Omega$ is finite and
    $\phi(Y)=[\DD{i_1}{j_1},\dotsc,\DD{i_k}{j_k}]$, we can
    partition~$Y$ as $Y=Y_1\cup \dotsb \cup Y_k$ with $|Y_l|=i_l$ for
    each~$l$.  Furthermore, if $Z\subseteq Y$ and $Z=Z_1\cup \dotsb
    \cup Z_k$, where $Z_l\subseteq Y_l$ for each~$l$, then we can
    write~$\phi(Z)$ as a sum of multisets $\phi(Z)=M_1+\dotsb+M_k$,
    where $\wt(M_l)=|Z_l|$ for each~$l$ and $M_l=[\DD{i_l}{j_l}]$ if
    $Z_l=Y_l$.
  \end{enumerate}
\end{defn}

Note that condition~(ii) implies that $\phi(Y_l)=[\DD{i_l}{j_l}]$ for
$j=1$, $2$, \dots,~$k$.  Essentially, this condition means that
subsets of~$Y\nomcorr$ correspond to ``submultisets'' of~$\phi(Y)$ in
a suitable sense.

In the case of $G=H\Wr S$ mentioned above, we take the connected blocks of
weight~$n$ to be the orbits of the action of~$H$ on $n$-subsets
of~$\Delta$.  Then every orbit of~$G$ can be put into correspondence
with a multiset of $H$-orbits as follows.  If $Y\subset\Omega$ is an
orbit representative, then
$\phi(Y)=[\,\pi_i(Y):\pi_i(Y)\ne\emptyset\,]$, where the $\pi_i$~are
projections: $\pi_i(Y)=\{\,\delta:(\delta,i)\in Y\,\}$, and we
identify orbits of~$H$ with orbit representatives.  Note that
$\wt(\phi(Y))=|Y|$ as required, and that condition~(ii) is also
satisfied; in fact, in the notation of the condition, we have
$M_l=[\DD{i'_l}{j'_l}]$ for each~$l$, for some appropriate $i'_l$
and~$j'_l$.

Another example is the automorphism group of the random graph.  The
random graph is the unique countable homogeneous structure whose age
consists of all finite graphs.  It is also known as the Fra{\"\i}ss\'e
limit of the set of finite graphs; see Cameron~\cite{CameronBook} for
more information on homogeneous structures and Fra{\"\i}ss\'e's
theorem.  We take the set of connected blocks to be the isomorphism
classes of finite connected graphs, where the weight of a connected
block is the number of vertices in it.  Any orbit can be uniquely
described by the multiset of connected graph components in an orbit
representative.  Condition~(i) is immediate, as is condition~(ii).
Note, however, that there are examples in this scenario where
$M_l$~may not be a singleton.  For example, if $Y=P_2$ is the path of
length~$2$ (with three vertices), so that $\phi(Y)=[P_2]$, and
$Z\subset Y$ consists of the two end vertices of the path, then
$\phi(Z)=[K_1,K_1]$.

This prototypical example explains the choice of terminology: the
basic units in this example are the connected graphs, so we have
called our basic units connected blocks, both to suggest this example
and that of strongly connected components in tournaments as considered
in section~\ref{sec:special-II} below.

Cameron~\cite[Sec.~2]{CameronModelThy} has shown that $A(G)$~is a
polynomial algebra if $G$~is an oligomorphic wreath-$S$-like group,
from which it follows that $A(G)$~is an integral domain in this case.
It also follows that $\varepsilon$~is a prime element, so both
Conjectures~\ref{conj:eprime} and~\ref{conj:intdomain} hold in this
case.  The argument that $A(G)$~is a polynomial algebra in the
oligomorphic case is similar to that presented below for
wreath-$A$-like groups, only significantly simpler.

We now show, using a new argument based on Ramsey-orderings, that
$A(G)$~is an integral domain in the wreath-$S$-like case, even without
the assumption that $G$~is oligomorphic.  This will also provide a
basis for the arguments presented in the next section for
wreath-$A$-like groups.

\begin{theorem}
  \label{thm:wreath-S}
  If $G\nomcorr$ is wreath-$S\nomcorr$-like, then $A(G)\!$ is an
  integral domain.
\end{theorem}

\begin{proof}
  We claim that in such a situation, the conditions of
  Conjecture~\ref{conj:main} are satisfied, and hence $A(G)$~is an
  integral domain.
  
  Following the requirements of the conjecture, let $m$ and~$n$ be
  positive integers and pick any integer $N\ge 2(m+n)$.  Denote the
  inverse of $\phi$ by $\psi$ and let $\alpha$ run through all
  multisets of connected blocks of total weight $m+n$, then we set
  $c_\alpha=\psi(\alpha)$ and let $X_\alpha$~be an $N$-set in the
  orbit $\psi(\alpha+[\DD11,\dotsc,\DD11])$, where the second multiset
  has $N-(m+n)$ copies of~$\DD11$.  We claim that this gives a Ramsey
  ordering of the orbits on $(m+n)$-sets, where the multisets are
  ordered lexicographically (which gives a well-ordering on the
  multisets).  Firstly, every $(m+n)$-set orbit appears among the list
  by hypothesis, as $\psi$~is a bijection.  Secondly, by construction,
  there is an $(m+n)$-subset of~$X_\alpha$ in the
  orbit~$\psi(\alpha)$, namely partition~$X_\alpha$ as in
  condition~(ii) of the definition, and remove all of the elements
  corresponding to the copies of~$\DD11$ added.  This subset will then
  map to~$\alpha$ under~$\phi$, by condition~(ii).  Finally, any
  $(m+n)$-subset of $X_\alpha$ can be seen to correspond to a multiset
  lexigraphically less than or equal to~$\alpha$, again using
  condition~(ii) and the fact that $\DD11$~is the least connected
  block, so the subset will be in an orbit~$c_\beta$ with
  $\beta\lexle\alpha$, as required.
  
  We note that the induced Ramsey orderings on $m$-set orbits and
  $n$-set orbits are given by precisely the same construction.
  Specifically, let $\beta$~be a multiset with $\wt(\beta)=n$.  Then
  the orbit corresponding to the multiset~$\beta$ first appears
  in~$X_{\alpha_0}$ where $\alpha_0=\beta+[\DD11,\dotsc,\DD11]$.  For
  assume that an $n$-set~$Z$ in the orbit~$\psi(\beta)$ appears in
  $X_\alpha$.  As we have $\phi(Z)=\beta$, $\beta$~must be a
  ``submultiset'' of~$\alpha$ in the sense of condition~(ii), and it
  is clear that the lexicographically smallest such~$\alpha$ is the
  one given by adjoining an appropriate number of copies of~$\DD11$
  to~$\beta$.  It is not difficult to show that $\beta\vee\gamma$ is
  precisely the multiset $\beta+\gamma$, and that $\beta\lexlt\beta'$
  implies $\beta+\gamma\lexlt\beta'+\gamma$, and therefore
  $\beta\vee\gamma\lexlt\beta'\vee\gamma$; similarly,
  $\gamma\lexlt\gamma'$ implies
  $\beta\vee\gamma\lexlt\beta\vee\gamma'$.  (The argument is similar
  to that of Theorem~\ref{thm:wreath-A} below.)  Thus the conditions
  of the conjecture are satisfied by this Ramsey ordering, and hence
  $A(G)$~is an integral domain.
\end{proof}

\ifthesis\clearpage\fi
\section{Special cases (II): Wreath-$A$-like groups}
\label{sec:special-II}

We can now apply the same ideas used for the wreath-$S$-like case to
the next class of groups, although the details are more intricate.
The only essential difference between these two classes is that here
we deal with ordered sequences of connected blocks instead of
unordered multisets of connected blocks.  We first define this class
of groups and show that their algebras are integral domains.  We then
show that in the oligomorphic case, they have a structure similar to
that of shuffle algebras, and deduce that they are polynomial rings.
With this information, we then look at some integer sequences which
arise from this family of groups.

\subsection{Wreath-$A$-like groups}
\label{sec:wreath-A}

If we have two finite sequences $S_1=(x_1,\dotsc,x_r)$ and
$S_2=(y_1,\dotsc,y_s)$, then we write $S_1\oplus
S_2=(x_1,\dotsc,x_r,y_1,\dotsc,y_s)$ for their concatenation.

\begin{defn}
  \label{defn:wr-A-like}
  We say that a permutation group~$G$ on~$\Omega$ is
  \emph{wreath-$A$-like} if there is a set of connected blocks~$\{\DD
  ij\}$ and a bijection~$\phi$ from the set of orbits of~$G$ on finite
  subsets of~$\Omega$ to the set of all finite sequences of connected
  blocks, with the bijection satisfying the following conditions:
  \begin{enumerate}
  \item  If $Y\subset\Omega$ is finite, then $\wt(\phi(Y))=|Y|$.
  \item If $Y\subset\Omega$ is finite and
    $\phi(Y)=(\DD{i_1}{j_1},\dotsc,\DD{i_k}{j_k})$, we can
    partition~$Y$ as an ordered union $Y=Y_1\cup \dotsb \cup Y_k$ with
    $|Y_l|=i_l$ for each~$l$.  Furthermore, if $Z\subseteq Y$ and
    $Z=Z_1\cup \dotsb \cup Z_k$, where $Z_l\subseteq Y_l$ for
    each~$l$, then we can write~$\phi(Z)$ as a concatenation of
    sequences $\phi(Z)=S_1\oplus\dotsb\oplus S_k$ where
    $\wt(S_l)=|Z_l|$ for each~$l$, and $S_l=(\DD{i_l}{j_l})$ if
    $Z_l=Y_l$.
  \end{enumerate}
\end{defn}

As in the wreath-$S$-like case, condition~(ii) implies that
$\phi(Y_l)=(\DD{i_l}{j_l})$ for $l=1$,~$2$, \dots,~$k$.

Our prototypical family of groups for this class of groups are those
of the form $G=H\Wr A$, where $H$~is a permutation group on~$\Delta$,
and $A$~is the group of all order-preserving permutations of the
rationals.  Again, the wreath product action is the imprimitive one,
so $G$~acts on $\Omega=\Delta\times\mathbb{Q}$.  As before, we take
the connected blocks of weight~$n$ to be the orbits of the action
of~$H$ on $n$-subsets of~$\Delta$.  Then every orbit of~$G$ can be put
into correspondence with a unique sequence of $H$-orbits as follows.
If $Y\subset\Omega$ is an orbit representative, we can apply an
element of the top group~$A$ to permute~$Y$ to a set of the form
$(\Delta_1\times\{1\})\cup(\Delta_2\times\{2\})\cup \dotsb
\cup(\Delta_t\times\{t\})$, where each $\Delta_i$~is non-empty.  Each
of the~$\Delta_i$ is a representative of some $H$-orbit, so we set
$\phi(Y)=(\Delta_1,\Delta_2,\dotsc,\Delta_t)$, again blurring the
distinction between orbits and orbit representatives.  It is again
easy to see that conditions (i) and~(ii) of the definition hold in
this case.

Another example is the automorphism group of the random tournament.
In this context, a tournament is a complete graph, every one of whose
edges is directed, and the random tournament is the Fra{\"\i}ss\'e
limit of the set of finite tournaments.  A tournament is called
strongly connected if there is a path between every ordered pair of
vertices.  It can be shown quite easily that every tournament can be
decomposed uniquely as a sequence of strongly connected components,
where the edges between components are all from earlier components to
later ones.  So here we take our set of connected blocks to be the
isomorphism classes of finite strongly connected tournaments (and
again, the weight of a connected block is the number of vertices in
it), and if $T$~is a finite subset of the random tournament, we set
$\phi(T)$~to be the sequence of strongly connected components
of~$T\nomcorr$.  Again, it is not difficult to see that conditions (i)
and~(ii) hold.  Also, as in the case of the random graph, it may be
that a sub-tournament has more components that the original
tournament; for example, the cyclically-oriented $3$-cycle is strongly
connected, but any $2$-element subset of it consists of two strongly
connected $1$-sets.

A third example is the automorphism group of the ``generic pair of
total orders''\negthinspace.  This is the Fra{\"\i}ss\'e limit of the
class of finite sets, where each finite set carries two (unrelated)
total orders, which can be taken as $a_1<a_2<\dotsb<a_n$ and
$a_{\pi(1)}\prec a_{\pi(2)}\prec\dotsb\prec a_{\pi(n)}$ for some
permutation $\pi\in S_n$.  Thus orbits of the Fra{\"\i}ss\'e limit are
described by permutations.  We can take the connected blocks for this
group to be the permutations $\pi\in S_n$ for which there exists
no~$k$ with $0<k<n$ such that $\pi$~maps $\{1,\dotsc,k\}$ to itself.
The details of this example are not hard to check.

\goodbreak
\begin{theorem}
  \label{thm:wreath-A}
  If $G\nomcorr$ is wreath-$A$-like, then $A(G)\!$ is an integral
  domain.
\end{theorem}

\begin{proof}
  The proof runs along very similar lines to that of
  Theorem~\ref{thm:wreath-S}.  If $\alpha$~is a sequence of connected
  blocks, we write~$[\alpha]$ to denote the multiset whose elements
  are the terms of the sequence with their multiplicities.  We define
  an ordering on sequences by $\alpha<\beta$ if
  $[\alpha]\lexlt[\beta]$ or $[\alpha]=[\beta]$ and
  $\alpha\lexgt\beta$.
  
  Again, we show that the conditions of Conjecture~\ref{conj:main} are
  satisfied in this case.  Let $m$ and~$n$ be positive integers and
  let $N$~be a positive integer with $N\ge 2(m+n)$.  Denoting the
  inverse of $\phi$ by $\psi$ and letting $\alpha$ run through all
  sequences of connected blocks of total weight $m+n$, we set
  $c_\alpha=\psi(\alpha)$ and let $X_\alpha$~be a $N$-set in the orbit
  $\psi(\alpha\oplus(\DD11,\dotsc,\DD11))$, where the second sequence
  has $N-(m+n)$ copies of~$\DD11$.  We claim that this gives a Ramsey
  ordering of the orbits on $(m+n)$-sets, where the sequences are
  ordered as described in the previous paragraph.  Firstly, every
  $(m+n)$-set orbit appears in the list by hypothesis, as $\psi$~is a
  bijection.  Secondly, by construction, there is an $(m+n)$-subset
  of~$X_\alpha$ in the orbit~$\psi(\alpha)$, namely
  partition~$X_\alpha$ as in condition~(ii) of the definition, and
  remove all of the elements corresponding to the copies of~$\DD11$
  appended.  This subset will then map to~$\alpha$ under~$\phi$, by
  condition~(ii).
  
  To show the final condition of Ramsey orderings, we must show that
  any $(m+n)$-subset of $X_\alpha$ is in an orbit corresponding to a
  sequence less than or equal to~$\alpha$.  Using the notation of
  condition~(ii), we let $\alpha=(\DD{i_1}{j_1},\dotsc,\DD{i_k}{j_k})$
  and $X_\alpha=X_1\cup \dotsb \cup X_k \cup X_{k+1} \cup \dotsb \cup
  X_r$, where $X_{k+1}$, \dots, $X_r$ correspond to the appended
  copies of~$\DD11$.  Consider a subset $Y=Y_1\cup\dotsb\cup
  Y_r\subset X_\alpha$ with $|Y|=m+n$.  If $Y_l\ne X_l$ for some~$l$
  with $X_l\ne\DD11$, then clearly $[\phi(Y)]\lexlt[\alpha]$, as
  $\wt(S_l)<i_l$, and the only new connected blocks which can be used
  are copies of~$\DD11$, which is the least connected block.  So the
  remaining case to consider is where some of the~$\DD{i_l}{j_l}$ are
  equal to~$\DD11$, and for some or all of those, $Y_l=\emptyset$,
  whereas $Y_s=X_s$ for some $s>k$.  But in such a case, while we have
  $[\phi(Y)]=[\alpha]$, it is clear that $\phi(Y)\lexge\alpha$.  So in
  either case, we have $\phi(Y)\le\alpha$, or equivalently $Y\le
  c_\alpha$, as required.

  We note that the induced Ramsey orderings on $m$-set orbits and
  $n$-set orbits are given by precisely the same construction; in
  particular, the orbit given by the sequence~$\beta$ first appears
  in~$X_\alpha$, where $\alpha=\beta\oplus(\DD11,\dotsc,\DD11)$.
  
  Finally, we must show that the remaining conditions of the
  conjecture are satisfied by this Ramsey ordering.  We will only show
  that $\beta<\beta'$ implies $\beta\vee\gamma<\beta'\vee\gamma$; the
  other condition follows identically.  We first deduce an explicit
  description of $\beta\vee\gamma$.
  
  A \textit{shuffle} of two sequences, say $(x_1,\dotsc,x_r)$ and
  $(y_1,\dotsc,y_s)$, is a sequence $(z_1,\dotsc,z_{r+s})$ for which
  there is a partition of $\{1,2,\dotsc,r+s\}$ into two disjoint
  sequences $1\le i_1<i_2<\dotsb<i_r\le r+s$ and $1\le
  j_1<j_2<\dotsb<j_s\le r+s$ with $z_{i_k}=x_k$ for $1\le k\le r$ and
  $z_{j_k}=y_k$ for $1\le k\le s$.
  
  We first show that $\beta\vee\gamma$ is the lexicographically
  greatest shuffle of~$\beta$ with~$\gamma$; this is not difficult
  although the argument is a little intricate.  We let $\alpha_0$~be
  this greatest shuffle and note that $[\alpha_0]=[\beta]+[\gamma]$.
  Now let $\alpha$~be any sequence of connected blocks for which
  $c_\alpha$ contains a $d_\beta\cup e_\gamma$ decomposition; we must
  show that $\alpha_0\le\alpha$.  (Here $d_\beta$~and~$e_\gamma$ are
  the orbits on $m$-sets and $n$-sets corresponding to $\beta$
  and~$\gamma$ respectively.)

  \newcommand{\pp}{^{\vphantom{\prime}}}
  We let $\alpha=(A_1,\dotsc,A_k)$ be this sequence of connected
  blocks, and let $Y$~be a representative of the orbit~$c_\alpha$.
  Write~$Y$ as an ordered union $Y=Y_1\cup\dotsb\cup Y_k$ as in
  condition~(ii) of the definition of wreath-$A$-like groups.  Then
  any decomposition of~$c_\alpha$ into two subsets can be written as
  $$c_\alpha=Z\cup Z'=(Z\pp_1\cup\dotsb\cup
  Z\pp_k)\cup(Z'_1\cup\dotsb\cup Z'_k),$$
  where $Y\pp_l=Z\pp_l\cup
  Z'_l$ as a disjoint union for each~$l$.  Now if we require
  $\phi(Z)=\beta$ and $\phi(Z')=\gamma$, this means that the sequences
  $S_1\oplus\dotsb\oplus S_k$ and $S'_1\oplus\dotsb\oplus S'_k$
  corresponding to~$Z$ and~$Z'$ respectively, as given by
  condition~(ii), must equal $\beta$ and~$\gamma$ respectively.  If
  $\{Z\pp_l,Z'_l\}=\{Y\pp_l,\emptyset\}$, then
  $[S\pp_l]+[S'_l]=[A\pp_l]$ by condition~(ii), but if not, then
  $[S\pp_l]+[S'_l]\lexlt[A\pp_l]$ by comparing weights.  As $M_1\lexlt
  M_2$ implies $M_1+M\lexlt M_2+M$ for any multisets $M_1$, $M_2$
  and~$M$, it follows that $[\beta]+[\gamma]\lexle[\alpha]$ with
  equality if and only if $\{Z\pp_l,Z'_l\}=\{Y'_l,\emptyset\}$ for
  each~$l$, that is, $[\alpha_0]\lexle[\alpha]$ with equality if and
  only if $\alpha$~is a shuffle of $\beta$ and~$\gamma$.  And if
  $\alpha$~is such a shuffle, then $\alpha\lexle\alpha_0$ by
  construction, so $\alpha_0\le\alpha$, as required.
  
  Given this, we can now show that if $\beta<\beta'$, then
  $\beta\vee\gamma<\beta'\vee\gamma$.  We first consider the case that
  $[\beta]\lexlt[\beta']$, from which it follows that
  $[\beta]+[\gamma]\lexlt[\beta']+[\gamma]$.  Since
  $[\beta\vee\gamma]=[\beta]+[\gamma]$ and
  $[\beta'\vee\gamma]=[\beta']+[\gamma]$, we deduce that
  $[\beta\vee\gamma]\lexlt[\beta'\vee\gamma]$, so
  $\beta\vee\gamma<\beta'\vee\gamma$.

  Now consider the other possible case, namely $[\beta]=[\beta']$ but
  $\beta\lexgt\beta'$.  Note that
  $[\beta\vee\gamma]=[\beta'\vee\gamma]$ in this case, so we must show
  that $\beta\vee\gamma\lexgt\beta'\vee\gamma$.  We let
  $\beta=(\Delta_1,\dotsc,\Delta_r)$,
  $\beta'=(\Delta'_1,\dotsc,\Delta'_r)$ and $\gamma=(E_1,\dotsc,E_s)$
  in the following.  We also let
  $\alpha=\beta\vee\gamma=(A_1,\dots,A_{r+s})$ and
  $\alpha'=\beta'\vee\gamma=(A'_1,\dotsc,A'_{r+s})$.  Recalling that
  $\beta\vee\gamma$ is the lexicographically greatest shuffle of
  $\beta$ and~$\gamma$, we can construct $\beta\vee\gamma$ by using
  the following merge-sort algorithm (written in pseudo-code).
  \newcommand{\from}{\leftarrow}
  \begin{tabbing}
    \quad\=\quad\=\quad\=\quad\=\kill
    \+\kill
    \textbf{function} \textit{MergeSort}$(\beta,\gamma)$\\
    \+\kill
    \>\>$\{$ \textit{We have $\beta=(\Delta_1,\dotsc,\Delta_r)$ and
      $\gamma=(E_1,\dotsc,E_s)$} $\}$\\
    $i\from 1$\\
    $j\from 1$\\
    \textbf{while} $i\le r$ \textbf{or} $j\le s$ \textbf{do}\\
    \>\textbf{if} $(i>r)$ \textbf{then}
      $\{$ $A_{i+j-1}\from E_j$; \ $j\from j+1$ $\}$\\
    \>\textbf{else if} $(j>s)$ \textbf{then}
      $\{$ $A_{i+j-1}\from\Delta_i$; \ $i\from i+1$ $\}$\\
    \>\textbf{else if} $(E_j\ge\Delta_i)$ \textbf{then}
      $\{$ $A_{i+j-1}\from E_j$; \ $j\from j+1$ $\}$\\
    \>\textbf{else}
      $\{$ $A_{i+j-1}\from\Delta_i$; \ $i\from i+1$ $\}$\\
    \textbf{od}\\
    \textbf{return} $\alpha=(A_1,\dotsc,A_{r+s})$
  \end{tabbing}
  Observe what happens if we run the algorithm on the pairs
  $(\beta,\gamma)$ and $(\beta',\gamma)$.  Assume that
  $\Delta_i=\Delta'_i$ for $i<i_0$, but that
  $\Delta_{i_0}>\Delta'_{i_0}$.  Then they will run identically as
  long as $i<i_0$.  When $i=i_0$, they will both continue taking terms
  from~$\gamma$ until $E_j<\Delta_{i_0}$ or $\gamma$~is exhausted.
  Once this happens, the $(\beta,\gamma)$ algorithm will
  take~$\Delta_{i_0}$ next, so $A_{i_0+j-1}=\Delta_{i_0}$, but the
  $(\beta',\gamma)$ algorithm will take $\max\{\Delta'_{i_0},E_j\}$,
  so
  $A'_{i_0+j-1}=\max\{\Delta'_{i_0},E_j\}<\Delta_{i_0}=A_{i_0+j-1}$.
  Thus we have $\beta\vee\gamma\lexgt\beta'\vee\gamma$, so
  $\beta\vee\gamma<\beta'\vee\gamma$ as required.

  It follows that $A(G)$~is an integral domain, as we wanted.
\end{proof}

\subsection{Shuffle algebras}
\label{sec:shuffle}

In the oligomorphic case, we can do better: the algebra~$A(G)$ is
actually a polynomial algebra if $G$~is an oligomorphic
wreath-$A$-like group.  We show this by noting strong similarities
between our algebra and standard shuffle algebras, and using
well-known properties of shuffle algebras, in particular that the
Lyndon words form a polynomial basis for the shuffle algebra.

We start by briefly recalling the key facts we will need.  We take
these results from Reutenauer's book on free Lie algebras~\cite{Reut}.
The references to definitions, theorems and so forth are to his book.

Let $\alphabet$ be an alphabet.  Although Reutenauer sometimes assumes
the alphabet to be finite, it will be clear that all of the results we
use below work equally well in the infinite case: since words are
always of finite length and we only ever work with finitely many words
at once, we can always restrict attention to the finite subset
of~$\alphabet$ containing the letters in use.

We write $\alphabet^*$~for the set of words in the
alphabet~$\alphabet$.  We write $\KA$~for the $\field$-vector space
with basis~$\alphabet^*$\negthinspace.  If we use the concatenation
product (where the product of two words is just their concatenation),
then this is the ring of non-commuting polynomials over~$\alphabet$.
But there is another product that we can define on words, and by
extension on~$\KA$, called the \textit{shuffle product}.  This is
explained in section~1.4 of Reutenauer, and we now essentially quote
parts of it.

Let $w=a_1\!\dotsm a_n$ be a word of length~$n$
in~$\alphabet^*$\negthinspace, and let $I\subseteq\{1,\dotsc,n\}$.  We
denote by~$w|I$ the word $a_{i_1}\!\dotsm a_{i_k}$ if
$I=\{i_1<i_2<\dotsb<i_k\}$; in particular, $w|I$~is the empty word if
$I=\emptyset$.  (Such a word~$w|I$ called a \textit{subword} of~$w$.)
Note that when
$$\{1,\dotsc,n\}=\bigcup_{j=1}^p I_j,$$
then $w$~is determined by the $p$~words $w|I_j$ and the $p$~subsets
$I_j$.

Given two words $u_1$ and~$u_2$ of respective lengths
$n_1$ and~$n_2$, their \textit{shuffle product}, denoted by
$u_1\shuffle u_2$, is the polynomial
$$u_1\shuffle u_2=\sum w(I_1,I_2),$$
where the sum is taken over all pairs $(I_1,I_2)$ of disjoint subsets
of $\{1,\dotsc,n\}$ with $I_1\cup I_2=\{1,\dotsc,n\}$ and $|I_j|=n_j$
for $j=1$,~$2$, and where the word $w=w(I_1,I_2)$ is defined by
$w|I_j=u_j$ for $j=1$,~$2$.  Note that $u_1\shuffle u_2$ is a sum of
words of length~$n$, each with the same multiset of letters, and so is
a homogeneous polynomial of degree~$n$.  Note also that the empty
word, denoted by~$1$, is the identity for the shuffle product, that
the shuffle product is commutative and associative, and that it is
distributive with respect to addition.  Thus $\KA$~with the shuffle
product is a commutative, associative algebra, called the
\textit{shuffle algebra}.

Using the associative and distributive properties of the shuffle
product, we can also give an expression for the shuffle product of the
words $u_1$, \dots,~$u_p$, of respective lengths $n_1$, \dots,~$n_p$;
their shuffle product is the polynomial
$$u_1\shuffle\dotsb\shuffle u_p=\sum w(I_1,\dotsc,I_p),$$
where now the sum is taken over all $p$-tuples $(I_1,\dotsc,I_p)$ of
pairwise disjoint subsets of $\{1,\dotsc,n\}$ with $\bigcup_{i=1}^p
I_j=\{1,\dotsc,n\}$ and $|I_j|=n_j$ for each $j=1$, \dots,~$p$, and
where the word $w=w(I_1,\dotsc,I_p)$ is defined by $w|I_j=u_j$ for
each $j=1$, \dots,~$p$.

A word appearing in the shuffle product $u_1\shuffle \dotsb\shuffle
u_p$ is called a \textit{shuffle} of $u_1$, \dots,~$u_p$.  Note that
this is consistent with the definition of shuffle we used in the proof
of Theorem~\ref{thm:wreath-A} above.  As an example, if
$a,b,c\in\alphabet$, then $ab\shuffle ac=abac+2aabc+2aacb+acab$, and
$aabc$ and $acab$ are both shuffles of $ab$ and~$ac$.

The next definition we need is that of a Lyndon word.  Assume that our
alphabet~$\alphabet$ is totally ordered.  Then a \textit{Lyndon word}
in~$\alphabet^*$ is a non-empty word which is lexicographically
smaller than all of its nontrivial proper right factors; in other
words, $w$~is a Lyndon word if $w\ne1$ and if for each factorisation
$w=uv$ (concatenation product) with $u,v\ne1$, one has $w\lexlt v$.

An alternative categorisation of Lyndon words is as follows
(Corollary~7.7 in Reutenauer).  Given a word $w=a_1\dotsm a_n$ of
length~$n$, we can define the rotation operator~$\rho$ by
$\rho(w)=a_2\dotsm a_na_1$.  Then a word~$w$ of length $n\ge1$ is
Lyndon if and only if $w\lexlt\rho^k(w)$ for $k=1$, \dots, $n-1$,
which is to say that $w$~is primitive (it does not have the form
$w=u^r$ for some $r>1$) and that it is lexicographically smaller than
any rotation (cyclic permutation) of itself.  It follows that Lyndon
words are in bijective correspondence with primitive necklaces; see
\cite[Chap.~7]{Reut} for more information.

A key property of Lyndon words is that every word $w\in\alphabet^*$
can be written \textit{uniquely} as a decreasing product of Lyndon
words, so $w=l_1^{r_1}\!\dotsm l_k^{r_k}$, where
$l_1\lexgt\dotsb\lexgt l_k$ and $r_1,\dotsc,r_k\ge1$.  (This follows
from Theorem~5.1 and Corollary~4.4, and can also easily be proved
directly\textemdash see section~7.3.)

Finally, Theorem 6.1 states that the shuffle algebra~$\KA$ is a
polynomial algebra generated by the Lyndon words, and that for each
word~$w$, written as a decreasing product of Lyndon words
$w=l_1^{r_1}\dotsm l_k^{r_k}$ as in the previous paragraph, one has
\begin{equation}
  \label{eq:lyndon-shuffle}
  S(w) \stackrel{\text{def}}{=}
    \frac{1}{r_1!\dotsm r_k!} \, l_1^{\shuffle r_1} \shuffle\dotsb
    \shuffle l_k^{\shuffle r_k}
    = w+\sum_{\substack{[u]=[w]\\u\lexlt w}} \alpha_u u,
\end{equation}
for some non-negative integers~$\alpha_u$, where $l^{\shuffle r}$
means $l\shuffle \dotsb\shuffle l$ with $r$~terms in the product, and,
in this context, $[u]$~means the multiset of letters in the word~$u$.

Note that it is equation~\eqref{eq:lyndon-shuffle} which proves that
$\KA$~is a polynomial algebra: the set~$\alphabet^*$ is a
$\field$-vector space basis for~$\KA$, and given any finite
multiset~$M$ of elements of~$\alphabet$, the matrix relating the basis
elements $\{\,w:w\in\alphabet^*\ \text{and}\ [w]=M\,\}$ to
$\{\,S(w):w\in\alphabet^*\ \text{and}\ [w]=M\,\}$ is unitriangular
when the words are listed in lexicographic order, so that
$\{\,S(w):w\in\alphabet^*\,\}$ also forms a basis for~$\KA$.  This
argument is true whether $\alphabet$~is finite or infinite.

We can now apply this to our case of oligomorphic wreath-$A$-like
permutation groups.  Let $G$ acting on~$\Omega$ be such a group, as in
Definition~\ref{defn:wr-A-like} above.  We obviously take our
alphabet~$\alphabet$ to be the set of connected blocks of the action
(as given by the definition of wreath-$A$-like groups), so that
$\alphabet^*$~corresponds bijectively to the set of orbits of~$G$ on
finite subsets of~$\Omega$.  The alphabet~$\alphabet$ has the standard
ordering defined on connected blocks, and the set~$\alphabet^*$ can then be
ordered either by the lexicographic order (denoted~$\lexlt$) or by the
order we defined at the start of Theorem~\ref{thm:wreath-A}
(denoted~$<$).

Clearly $A(G)$ can be regarded as a $\field$-vector space, with the
set of characteristic functions of finite orbits as basis.  We will
identify the connected block sequence
$w=(\DD{i_1}{j_1},\dotsc,\DD{i_k}{j_k})$ with the characteristic
function of the corresponding orbit, writing $w$~for both.  Via this
correspondence, we can identify $A(G)$ with~$\KA$ as vector spaces.
The grading on~$A(G)$ induces a grading on~$\KA$: the homogeneous
component~$V_n(G)$ is identified with the subspace of~$\KA$ spanned by
$\{\,w\in\alphabet^*:\wt(w)=n\,\}$.  We then consider the product that
the vector space~$\KA$ inherits via this identification.  Let
$v\in\alphabet^*$ be another connected block sequence.  We write
$v\cshuffle w$ for the product in~$A(G)$ and the induced product
in~$\KA$.  The notation is designed to indicate that this product is
related to the shuffle product, as we will see, and we call it the
\textit{complete shuffle product}.  (It is also somewhat related to
the infiltration product on~$\KA$; see \cite[sect.~6.3]{Reut}.)
Recalling the definition of multiplication in~$A(G)$, we see that for
any finite subset $X\subset\Omega$ with $|X|=\wt(v)+\wt(w)$,
$$(v\cshuffle w)(X)=\sum_{\substack{Y\subseteq X\\|Y|=\wt(v)}}
v(Y)w(X\setminus Y).$$
But $v(Y)$ is none other than the
characteristic function which has value~$1$ if $\phi(Y)=v$ and
$0$~otherwise, and similarly for $w(Y\setminus X)$.  So we have
$$(v\cshuffle w)(X)=|\{\,Y\subseteq X:\phi(Y)=v,\ \phi(X\setminus
Y)=w\,\}|.$$
Thus, setting $u=\phi(X)$ and writing $u\to v\cup w$ if
there is a $Y\subseteq X$ with $\phi(Y)=v$ and $\phi(X\setminus Y)=w$,
we have
$$v\cshuffle w=\sum_{u\in\alphabet^*} \beta_u u,$$
where $\beta_u>0$ if $u\to v\cup w$ and $\beta_u=0$ otherwise.

Now we can characterise those $u$ for which $u\to w\cup v$ quite
easily.  Firstly, consider the case that $[u]=[w]+[v]$, that is, the
set of connected blocks of~$u$ is the same as those of $w$ and~$v$
combined.  Then $u\to w\cup v$ if and only if $u$~is a shuffle of $w$
and~$v$, by condition~(ii) of Definition~\ref{defn:wr-A-like}, as in
the proof of Theorem~\ref{thm:wreath-A}.  In fact, the terms
in~$w\cshuffle v$ with $[u]=[w]+[v]$ will be precisely $w\shuffle v$,
which is easy to see.  Now consider those terms with $[u]\ne[w]+[v]$.
If $[u]\lexlt[w]+[v]$, then it is easy to see that we cannot have
$u\to w\cup v$, but it may be possible otherwise.  We deduce that our
product is given by:
\begin{equation}
  \label{eq:complete-shuffle}
  w\cshuffle v = w \shuffle v +
  \sum_{\substack{\wt(u)=\wt(w)+\wt(v)\\ [u]\lexgt[w]+[v]}}
  \beta_u u
\end{equation}
for some non-negative integers~$\beta_u$.

Now given $w=l_1^{r_1}\!\dotsm l_k^{r_k}$ written as a (concatenation)
product of decreasing Lyndon words, we can consider the complete
shuffle product as we did for the normal shuffle product above:
\begin{equation}
  \label{eq:lyndon-cshuffle}
  \begin{split}
  \bar S(w)
  &\stackrel{\text{def}}{=} \frac{1}{r_1!\dotsm r_k!} \,
    l_1^{\cshuffle r_1} \cshuffle\dotsb \cshuffle l_k^{\cshuffle r_k} \\
  &\stackrel{\phantom{\text{def}}}{=} \frac{1}{r_1!\dotsm r_k!} \,
    l_1^{\shuffle r_1} \shuffle\dotsb \shuffle l_k^{\shuffle r_k} +
    \sum_{\substack{\wt(u)=\wt(w)\\ [u]\lexgt[w]}} \beta_u u\\
  &\stackrel{\phantom{\text{def}}}{=}
    w+\sum_{\substack{[u]=[w]\\u\lexlt w}} \alpha_u u +
    \sum_{\substack{\wt(u)=\wt(w)\\ [u]\lexgt[w]}} \beta_u u\\
  &\stackrel{\phantom{\text{def}}}{=}
    w+\sum_{\substack{\wt(u)=\wt(w)\\u>w}} \alpha_u u,
  \end{split}
\end{equation}
where the $\alpha_u$ and the $\beta_u$ are non-negative integers.  To
get the second line, we have repeatedly used
equation~\eqref{eq:complete-shuffle} to reduce the complete shuffle
product to a normal shuffle product.  Observe that
$\wt(l_1^{r_1}\dotsm l_k^{r_k})=\wt(w)$, hence the sum is over words
with $\wt(u)=\wt(w)$, and with $[u]\lexgt[w]$, since $\lexgt$~is
transitive and $[u_1]\lexgt[u_2]$ implies $[u_1]+[u]\lexgt[u_2]+[u]$
for any word~$u$.  That the $\beta_u$~are non-negative is easy to see,
and it is not that much harder to see that they are integral, although
we do not need this.  In the third line, we have used
equation~\eqref{eq:lyndon-shuffle}, and in the last line, we have set
$\alpha_u=\beta_u$ in the case that $[u]\lexgt[w]$, and used the
relation on words (sequences) defined in the previous section, namely
$u>w$ if $[u]\lexgt[w]$ or $[u]=[w]$ and $u\lexlt w$.

It is also important to note that in our case, the set
$\{\,u:\wt(u)=\wt(w)\,\}$ is finite, as there are only finitely many
connected blocks of each weight, the same number as the number of
orbits on sets of size $\wt(w)$, so that the sums in
equation~\eqref{eq:lyndon-cshuffle} are all finite.

We now see, as above, that the matrix relating $\{\,w:w\in\alphabet^*\ 
\text{and}\;\wt(w)=n\,\}$ to $\{\,\bar S(w):w\in\alphabet^*\ \text{and}\;
\wt(w)=n\,\}$ is unitriangular when the words of weight~$n$ are listed
in the order we have defined.  It follows that the~$\bar S(w)$ form a
vector space basis for~$A(G)=\KA$, and hence the set of Lyndon words
is a set of polynomial generators for~$A(G)$.  We summarise these
results as a theorem.

\vspace{0pt plus 1cm}\goodbreak

\begin{theorem}
  \label{thm:wreath-A-poly}
  If $G\nomcorr$ is an oligomorphic wreath-$A$-like permutation group,
  then $A(G)\!$ is a polynomial ring, and the generators are those
  characteristic functions on orbits corresponding to Lyndon words as
  described above.\qed
\end{theorem}

We can now deduce:

\begin{corollary}
  If $G\nomcorr$ is an oligomorphic wreath-$A$-like permutation group,
  then the element $\varepsilon\in V_1(G)\!$ is prime in~$A(G)$.
\end{corollary}

\begin{proof}
  We have $e=\DD11+\dotsb+\DD1r$, where the $\DD1j$~are the orbits on
  $1$-sets.  As each of the~$\DD1j$ is a Lyndon word,
  $A(G)=\field[\DD11,\dotsc,\DD1r,\DD21,\dotsc]$.  It follows that we
  can replace the polynomial generator~$\DD11$ by~$\varepsilon$ (as
  they are linearly related), giving
  $A(G)=\field[\varepsilon,\DD12,\dotsc,\DD1r,\DD21,\dotsc]$.  It is
  clear, since we then have
  $A(G)/(\varepsilon)\cong\field[\DD12,\dotsc,\DD1r,\DD21,\dotsc]$,
  that $A(G)/(\varepsilon)$ is an integral domain, so $\varepsilon$~is
  prime in~$A(G)$.
\end{proof}

\subsection{Integer sequences, necklaces and free Lie algebras}
\label{sec:sequences}

Theorem~\ref{thm:wreath-A-poly} leads us to revisit some counting
questions.  Cameron~\cite{CameronJIS} considered the following
question.  If the algebra~$A(G)$ corresponding to an ``interesting''
oligomorphic group~$G$ were polynomial, what would be the sequence
counting the number of polynomial generators of each degree?  From
knowledge of the dimension of each homogeneous component of~$A(G)$,
the answer can be determined using the inverse Euler transform.  Now
that we have an explicit description of the polynomial generators in
the wreath-$A$-like case, an examination of the sequences observed
might yield some interesting new information about those sequences.

The two sequences we will consider are those arising from the groups
$S_2\Wr A$ and $A\Wr A$, both of which appear in the On-Line
Encyclopedia of Integer Sequences~\cite{IntSeq}.  There are some
obvious generalisations to other groups, as we observe below.  The
$n$-th homogeneous component of the group $S_2\Wr A$ has
dimension~$F_{n+1}$ (a Fibonacci number, where $F_0=0$ and $F_1=1$),
and so the sequence counting the number of generators of degree~$n$ is
\htmladdnormallink{A006206}
{http://www.research.att.com/cgi-bin/access.cgi/as/njas/sequences/eisA.cgi?Anum=006206},
beginning $1$,~$1$, $1$, $1$, $2$, $4$, $5$, $8$, $11$, $18$,~\dots.
By our result, the $n$-th term of this sequence gives the number of
Lyndon words of weight~$n$ (starting with $n=1$) in the alphabet
$\alphabet=\{\Delta_1,\Delta_2\}$, where $\Delta_1$ and~$\Delta_2$ have
respective weights $1$ and~$2$.

Similarly, for the group $A\Wr A$, the $n$-th homogeneous component
has dimension~$2^{n-1}$ for $n\ge1$, and the sequence counting the
number of generators of degree~$n$ is
\htmladdnormallink{A059966}
  {http://www.research.att.com/cgi-bin/access.cgi/as/njas/sequences/eisA.cgi?Anum=059966},
beginning $1$,~$1$, $2$, $3$, $6$, $9$, $18$, $30$,~\dots.  (Note that
the paper quoted above had sequence
\htmladdnormallink{A001037}
  {http://www.research.att.com/cgi-bin/access.cgi/as/njas/sequences/eisA.cgi?Anum=001037}
by mistake, this being the inverse Euler transform of the closely
related sequence~$(2^n)$.)  This sequence then counts the number of
Lyndon words of weight~$n$ in the alphabet
$\alphabet=\{\Delta_1,\Delta_2,\dotsc\}$, where $\Delta_i$~has weight~$i$.

The Encyclopedia entry gives a different explanation, however: this
sequence lists the dimensions of the homogeneous components of the
free Lie algebra with one generator of each degree $1$, $2$, $3$,~etc.
The connection between these two descriptions of this sequence is easy
to describe, using \cite[Thm.\ 4.9]{Reut}.  Let $\alphabet$~be an
alphabet whose letters each have a positive integral degree\slash
weight (we use these terms interchangeably in this section), and where
there are only finitely many letters of each possible weight.  There
is a basis of the free Lie algebra on the alphabet~$\alphabet$ (viewed
as a vector space) given by $\{\,P_w:w\in\alphabet^*\ 
\text{Lyndon}\,\}$, where $P_a=a$ if $a\in\alphabet$, and
$P_w=[P_u,P_v]$ otherwise, where $w=uv$ with $v$~being the
lexicographically smallest nontrivial proper right factor of~$w$ (see
\cite[Thm.\ 5.1]{Reut}).  Note that it trivially follows by induction
that the degree of the homogeneous polynomial~$P_w$ is $\wt(w)$.  Thus
the dimension of the homogeneous component of degree~$n$ of the free
Lie algebra on the alphabet~$\alphabet$ is the number of Lyndon words
in~$\alphabet^*$ of weight~$n$.  It follows that we can also describe
the two sequences above as either the number of Lyndon words of
weight~$n$ in the alphabets $\{\Delta_1,\Delta_2\}$ and
$\{\Delta_1,\Delta_2,\dotsc\}$ respectively, or as the number of
primitive necklaces of weight~$n$ in these symbols, or as the
dimension of the homogeneous component of degree~$n$ of the free Lie
algebras on these sets.  This obviously generalises to other
wreath-$A$-like groups.

We may ask other counting questions based on these ideas.  We start
with an alphabet of weighted letters~$\alphabet$ (again with only
finitely many letters of each weight).  The primary questions arising
are how to transform between the three sequences:
\begin{align*}
  a_n &= \text{number of letters of weight $n$ in $\alphabet$,}\\
  w_n &= \text{number of words of weight $n$ in $\alphabet^*$,}\\
  l_n &= \text{number of Lyndon words of weight $n$ in $\alphabet^*\!$.}
\end{align*}
(Of course, $l_n$ can also be regarded as the number of primitive
necklaces of weight~$n$ in this alphabet.)  In our context, $a_n$ is
the number of connected blocks of weight~$n$ in our wreath-$A$-like
group, $w_n$ gives the dimension of the homogeneous component of
weight~$n$ in~$A(G)$ and $l_n$~gives the number of polynomial
generators of weight~$n$ in~$A(G)$.  We use the notation and some of
the ideas presented in Bernstein and Sloane's paper on integer
sequences~\cite{Mira}.

The transformation between $(a_n)$ and~$(w_n)$ can be effected by
INVERT, as every word is an ordered sequence of letters:
$$1+\sum_{n=1}^\infty w_nx^n = \frac{1}{1-\sum_{n=1}^\infty a_nx^n}.$$
The transformation between $(w_n)$ and~$(l_n)$ is performed using
EULER, as every word is a product of a decreasing sequence of Lyndon
words, so can be identified with a multiset of Lyndon words:
$$1+\sum_{n=1}^\infty w_nx^n = \prod_{n=1}^\infty
\frac{1}{(1-x^n)^{l_n}}.$$
It follows that we can transform between $(a_n)$ and $(l_n)$ using a
variant of WEIGH:
\begin{equation}
  \label{eq:count-lyndon}
  1-\sum_{n=1}^\infty a_nx^n = \prod_{n=1}^\infty (1-x^n)^{l_n}.
\end{equation}
Most of the six possible conversions between $(a_n)$, $(w_n)$
and~$(l_n)$ are straightforward given these formul\ae; the two which
are harder are converting $(w_n)$ and~$(a_n)$ to~$(l_n)$.  Inverting
the EULER transform is explained in~\cite{Mira}; we apply the same
idea to convert from~$(a_n)$ to~$(l_n)$.

Given a sequence $(a_n)$, we introduce the auxiliary sequence~$(c_n)$
defined by the equation $1-\sum_{n=1}^\infty a_nx^n =
\exp\bigl(-\sum_{n=1\mathstrut}^\infty c_nx^n\!/n\bigr)$.  Using the
generating functions $A(x)=\sum_{n=1}^\infty a_nx^n$ and
$C(x)=\sum_{n=1}^\infty c_nx^n$, we can perform standard manipulations
using the defining equation for~$(c_n)$ to deduce that
$C(x)=xA'(x)+C(x)A(x)$.  It follows that
\begin{equation}
  \label{eq:an->cn}
  c_n = na_n + \sum_{k=1}^{n-1} c_ka_{n-k}.
\end{equation}
Now substituting $\exp\bigl(-\sum c_nx^n/n\bigr)$ for $1-\sum a_nx^n$ in
equation~\eqref{eq:count-lyndon}, taking logarithms and expanding as a
power series gives the coefficient of $x^n\!/n$ to be
$c_n=\sum_{d\mid n} d\,l_d.$  Finally, M\"obius inversion gives
\begin{equation}
  \label{eq:cn->ln}
  l_n = \frac{1}{n}\,\sum_{d\mid n} \mu(n/d) c_d.
\end{equation}

Thus we have an effective way of calculating the number of Lyndon
words of a given weight given the number of letters of each possible
weight.

As an interesting example of this process, let us consider our
favourite group, $G=S_2\Wr A$.  In this case, recall that we have
$\alphabet=\{\Delta_1,\Delta_2\}$, so $a_1=a_2=1$ and $a_n=0$ for $n\ge3$.
Then the sequence~$(c_n)$ is calculated by equation~\eqref{eq:an->cn}:
we have $c_1=1$ and $c_2=3$.  For $n\ge3$, we have
$c_n=c_{n-1}+c_{n-2}$, so $(c_n)$~is the standard Lucas
sequence~$(L_n)$: $1$,~$3$, $4$, $7$, $11$, $18$,~\dots.  We can now
calculate the sequence~$(l_n)$: the first few terms are as we
predicted: $1$,~$1$, $1$, $1$, $2$, $2$, $4$, $5$,~\dots, and a
general formula is $l_n=\frac{1}{n}\sum_{d\mid n}\mu(n/d)L_d$, as is
given in the Encyclopedia entry for \htmladdnormallink{A006206}
{http://www.research.att.com/cgi-bin/access.cgi/as/njas/sequences/eisA.cgi?Anum=006206}.
One interesting thing to observe is that if $p$~is prime, then we
have $l_p=(\mu(1)L_p+\mu(p)L_1)/p=(L_p-1)/p$.  It follows that the
Lucas sequence satisfies $L_p\equiv1\pmod{p}$ for all primes~$p$, a
known result (see Hoggart and Bicknell~\cite{Hoggart}), but somewhat
surprising in this context.

The description of our sequence \htmladdnormallink{A006206}
{http://www.research.att.com/cgi-bin/access.cgi/as/njas/sequences/eisA.cgi?Anum=006206}
in the Encylcopedia is ``aperiodic binary necklaces [of
length~$n$] with no subsequence~$00$, excluding the sequence~`$0$'.''
Our description is that it counts primitive necklaces of weight~$n$ in
the alphabet $\{\Delta_1,\Delta_2\}$.  These are easily seen to be
equivalent: if we replace every~$\Delta_1$ by the symbol~$1$ and
every~$\Delta_2$ by the symbols~$10$ (in clockwise order, say), then
we will get a primitive (aperiodic) binary necklace with no
subsequence~$00$ whose length equals the weight of the necklace we
started with, and we can perform the inverse transformation equally
simply (as we are excluding the necklace~$0$).  We can do the same
with the group $S_n\Wr A$, enabling us to count the number of
primitive binary necklaces of length~$n$ with no subsequence
$00\dotsm0$ (with $n$~zeros) and excluding the necklace~$0$.

Now let us apply these ideas to the case $G=A\Wr A$.  Firstly, the
auxilary sequence turns out to be $c_n=2^n-1$, and the
sequence~$(l_n)$ is given by $l_n=\sum_{d\mid n}\mu(n/d)(2^d-1)$.
This can be simplified using the result $\sum_{d\mid n} \mu(n/d)=
[n=1]$, where we are using Iverson's convention that if $P$~is a
predicate, then $[P]=1$ if $P$~is true and $0$~otherwise.  So we have
$l_n=\sum_{d\mid n}\mu(n/d)2^d-[n=\nobreak1]$.  The sequence given by
$\sum_{d\mid n}\mu(n/d)2^d$ is sequence \htmladdnormallink{A001037}
{http://www.research.att.com/cgi-bin/access.cgi/as/njas/sequences/eisA.cgi?Anum=001037},
and so our sequence differs from it by~$1$ in the $n=1$ term only,
yielding the observed sequence \htmladdnormallink{A059966}
{http://www.research.att.com/cgi-bin/access.cgi/as/njas/sequences/eisA.cgi?Anum=059966}.
We can also give a necklace description of this sequence as above: it
is the number of primitive binary necklaces of length~$n$ excluding
the necklace~$0$---the sequence \htmladdnormallink{A001037}
{http://www.research.att.com/cgi-bin/access.cgi/as/njas/sequences/eisA.cgi?Anum=001037}
is essentially the same, but does not exclude the necklace~$0$, so it
it also counts the number of binary Lyndon words of
\textit{length}~$n$.  (These are the descriptions of this sequence
given in the Encyclopedia.)  Finally, as above, if we consider the
term~$l_p$ for $p$~prime, we see that
$l_p=((2^p-1)-1)/p=2(2^{p-1}-1)/p$, so for $p>2$, we deduce Fermat's
little theorem for base~$2$, that is $2^{p-1}\equiv1\pmod{p}$.

An investigation of those sequences of non-negative integers~$(b_n)$
for which $\frac{1}{n}\sum_{d\mid n}\mu(n/d)b_d$ is a non-negative
integer for all~$n$ has been undertaken by Puri and
Ward~\cite{PuriWard}, who call them exactly realizable.  We can thus
add to their work a class of exactly realizable sequences: those which
are of the form~$(c_n)$, where $(c_n)$~is given by
equation~\eqref{eq:an->cn} for some sequence of non-negative
integers~$(a_n)$.  A particular family of such sequences is given by
$a_i=1$ for $1\le i\le n$ and $a_i=0$ for $i>n$; these are sometimes
known as ``generalised Fibonacci sequences'', and have been discussed
by Du~\cite{Du} (where this sequence is called~$\phi_n$).  It would be
interesting to know whether new congruence identities can be
discovered by applying this technique to some of the sequences
identified there or to sequences produced by other wreath-$A$-like
groups.

\ifthesis\clearpage\fi
\section{Non-oligomorphic groups}
\label{sec:non-olig}

Throughout this \ifthesis part of the thesis\else paper\fi, we have
mostly focused on oligomorphic groups, proving results in general
where there was no problem in doing so.  In this final section, we
consider briefly the issues arising in the non-oligomorphic case.

As has already been pointed out above, the group~$\ZZ$ acting
regularly on~$\ZZ$ does not have a Ramsey ordering on $2$-sets, so
much of what we did above will not help us to understand the
algebra~$A(\ZZ)$.  It is easy to construct other similar examples.

A more difficult question is whether we have even got the ``right''
definition of the algebra~$A(G)$ in the non-oligomorphic case.  The
definition we have been using was introduced specifically to study the
behaviour of oligomorphic groups.  There are two finiteness conditions
which can be imposed on the algebra we consider.

Firstly, we have taken the direct sum $A(G)=\bigoplus_{n=0}^\infty
V_n(G)$, which is the direct limit as $N\to\infty$ of the vector
spaces $\bigoplus_{n=0}^N V_n(G)$ (with the obvious direct maps).  We
could have instead taken the cartesian sum $\sum_{n=0}^\infty V_n(G)$,
being the inverse limit of the same family of vector spaces (with the
obvious inverse maps).

Secondly, and independently of the first choice, we could either
take~$V_n(G)$ to be the vector space of all functions from $n$-subsets
of~$\Omega$ to~$\field$ which are fixed by~$G$, as we have until now,
or we could take it to be the subspace of this consisting of those
functions which assume only finitely many distinct values on $n$-sets.
(The latter idea was suggested to me by Peter Cameron.)  Note, though,
that if there are infinitely many orbits on $n$-sets, this vector
space will still have uncountable dimension.  It is not hard to check
that if we use the latter definition, the multiplication in the
algebra is still well-defined.  Also, this distinction does not exist
in the oligomorphic case.  (Another seemingly plausible choice, those
functions in~$V_n(G)$ which are non-zero on only finitely many orbits
of~$G$, can fail to produce a well-defined multiplication: consider,
for example, the case of~$e^2$ with our favourite non-oligomorphic
group,~$\ZZ$: it takes the value~$2$ on every $2$-set.)

Thus we have four plausible algebras to choose from, and it is not
clear which is the ``correct'' one to use.  More work is still
required in this area.

\ifthesis\else
\section*{Acknowledgements}

I would like to thank Prof.\ Peter Cameron, my supervisor, for
introducing me to this problem of his and giving me many, many helpful
pointers, hints and examples of groups which would fit into my
classes.  Also, thanks to Prof.\ Roger Bryant for his detailed
comments which helped to improve the exposition.
\fi

\ifthesis\clearpage\fi
\bibliographystyle{amsplain}
\bibliography{pgrefs}


\let\clearpage=\saveclearpage


\end{document}